\documentclass[11pt]{amsart}
\usepackage{amsmath,amssymb,enumerate}
\usepackage{epsf,epsfig,amsfonts,graphicx,color}
 \usepackage[applemac]{inputenc} 
\usepackage{enumerate}

 \date{}   

\theoremstyle{definition}

\newtheorem{remark}{\rm\bf Remark}

\newcommand{\weg}[1]{}
\newcommand{\be}{\begin{equation}}
\newcommand{\ee}{\end{equation}}

\title{A counterexample to Belgun-Moroianu conjecture} \date{}
\author{Vladimir S. Matveev   and Yuri Nikolayevsky }

\begin{document}
\address{Institute of Mathematics,  Friedrich-Schiller-Universit\"at Jena,  07737 Jena, Germany}
\address{Department of Mathematics and Statistics, La Trobe University, Melbourne, 3086, VIC, Australia}

\begin{abstract} We construct  an example of  a closed manifold with
a nonflat reducible locally metric  connection  such that it  preserves a conformal structure and such that it  is not the Levi-Civita connection of a Riemannian metric.
\end{abstract}

\maketitle

 \subsection*{ Definitions and result.}
All objects in our paper are assumed to be sufficiently smooth. We say that an affine torsion-free connection $\nabla$ on $M^n$ is \emph{locally metric}, if in a neighborhood of any $p \in M^n$ there exists a Riemannian metric $g$ which is parallel with respect to $g$.
We say that the connection \emph{preserves a conformal structure}, if there exists a Riemannian metric $g'$ on $M^n$ such that  for any vector field $V$ there exists a function $\mu$ such that $\nabla_Vg'= \mu g'$. We say that the connection is reducible, if its holonomy group $H_p\subseteq \mathrm{GL}(T_pM)$ is reducible, i.e., has  a nontrivial invariant subspace.

\vspace{1ex}

In \cite[Conjecture 1.3]{BM} it was conjectured by F. Belgun and  A. Moroianu
that \emph{on a closed manifold $M^n $ any reducible locally metric connection  that preserves a conformal structure is the Levi-Civita connection of a certain Riemannian metric, or is flat}.
 The main result of this paper is a counterexample to this statement.

\subsection*{\bf A counterexample to the Belgun-Moroianu conjecture.}
We take any  matrix $A= \begin{pmatrix}a_{11} & a_{12} \\ a_{21} & a_{22} \end{pmatrix} \in \mathrm{SL}(2, \mathbb{Z})$ such that  its eigenvalues are real positive and one of the eigenvalues which we denote by $\lambda$ is greater than $1$. The other eigenvalue is clearly $\tfrac{1}{\lambda}$.  For example, one can take $A= \begin{pmatrix}2 & 1 \\ 1& 1 \end{pmatrix}$.

We take the  2-torus $T^2= \mathbb{R}^2/{\mathbb{Z}^2}$ with the standard coordinates $x \textrm{ mod } 1$,  $y \textrm{ mod } 1$ and  $\mathbb{R}_+:= \{z\in \mathbb{R} \mid z >0\}$.  Consider the action of the group $\mathbb{Z}$ on the direct product $T^2 \times \mathbb{R}_+$ generated by
\begin{equation} \label{ac3}
f: T^2 \times \mathbb{R}_+ \to T^2 \times \mathbb{R}_+\ ,\  \
\begin{pmatrix} x\\ y\\ z\end{pmatrix} \stackrel{f}{\longmapsto} \begin{pmatrix}  A \left(\begin{array}{c} x\\ y\end{array}\right)  \\ \lambda  z\end{pmatrix}= \begin{pmatrix} a_{11} x + a_{12} y\\ a_{21} x + a_{22} y\\ \lambda z\end{pmatrix}.
\end{equation}
As $A \in \mathrm{SL}(2,\mathbb{Z})$, the mapping $f$ is well-defined and invertible.

Since this action of $\mathbb{Z}$ is discrete and free, the quotient $M^3:=  \left(T^2 \times \mathbb{R}_+\right)/{\mathbb{Z}}$ is a manifold. It is clearly closed (=compact)  since it is  homeomorphic to the following compact topological space: take the product of the torus  $T^2$ and
of  the interval $[1,\lambda]$, and identify each  point  $\begin{pmatrix}x\\ y\\1 \end{pmatrix}$ of   $T^2 \times \{1\}$ with the point $\begin{pmatrix}  A \left(\begin{array}{c} x\\ y\end{array}\right)  \\ \lambda  \end{pmatrix}$ of  $T^2 \times  \{\lambda\}$.

On $M^3$ we now construct a locally metric reducible nonflat connection.

Denote by $v_1$ an eigenvector of $A$ corresponding to the eigenvalue $\lambda$ and by $v_2$ an eigenvector corresponding to $\tfrac{1}{\lambda}$. We think of $v_1, v_2$ as of vector fields on the whole $T^2\times \mathbb{R}_+$ by extending them by translations (so that $v_1, v_2$ have constant components relative to the frame $\tfrac{\partial }{\partial x}, \tfrac{\partial }{\partial y}, \tfrac{\partial }{\partial z}$).  
Denote by $v_3$ the vector $\tfrac{\partial }{\partial z}$.

Consider the metric $g$ on $T^2\times \mathbb{R}_+$ such that at the point $(x,y,z)$ relative to the basis $v_1,v_2, v_3$ it is given by the matrix
$$
\begin{pmatrix} 1 & 0 & 0 \\ 0 & z^4 & 0 \\ 0 & 0 & 1\end{pmatrix}.
$$
An easy calculation shows that $f$ from \eqref{ac3} acts by the homothecy of $g$ with the coefficient $\lambda$: indeed,
\begin{equation*}
df(v_1)= Av_1= \lambda v_1, \; df(v_2)= Av_2=\tfrac{1}{\lambda} v_2  \text{ and } df(v_3)= \lambda v_3,
\end{equation*}
which shows that the push-forwards $df(v_i)$ of the vectors $v_1, v_2, v_3$ are mutually orthogonal, and also that the $g$-length of the push-forward of $v_i$ is $\lambda$ times the $g$-length of $v_i$.

Then both the Levi-Civita connection of $g$ and the conformal structure defined by $g$ are invariant with respect to the action of $f$. Therefore the former induces a connection $\nabla$ on $M^3$ which is a locally metric connection preserving the conformal structure.

To see that the connection $\nabla$ is reducible, we observe that the vector fields $v_1,v_2, v_3$ commute and therefore induce a (local) coordinate system $(\tilde x, \tilde y, z)$ on $T^2\times \mathbb{R}_+$, where the fact that the third coordinate is just $z$ follows from $v_3= \tfrac{\partial }{\partial z}$. Relative to these coordinates, the metric $g$  has  the form 
\begin{equation*} 
g= d\tilde x^2 + z^4 d \tilde y^2 + dz^2
\end{equation*}
and we see that the vector field $v_1= \tfrac{\partial }{\partial \tilde x}$ is parallel. Since the action of $\mathbb{Z}$ preserves  the 1-dimensional distribution spanned by $v_1$, this 1-dimensional distribution is parallel and is therefore invariant with respect to the holonomy group. A direct calculation shows that the metric $g$ and therefore the connection $\nabla$ are not flat. Since the metric $g$ on $T^2\times \mathbb{R}_+$ is not complete, the connection $\nabla$ is also not complete and hence can not be the Levi-Civita connection of a Riemannian metric on $M^3$.

Thus, we constructed an example of a closed manifold with a nonflat reducible torsion-free affine connection which is locally metric and preserves the conformal structure, but is not the Levi-Civita connection of any Riemannian metric.

\begin{remark} It is easy to generalize the example for any dimension $n\ge 3$. In dimension $n=2$, the Belgun-Moroianu conjecture is clearly true since any 2-dimensional locally metric reducible connection is flat.
\end{remark}

\begin{remark} The ideas  used in the construction of our counterexample came from the theory of compact solvmanifolds.
\end{remark}

\begin{remark} The universal covering space of $M^3$ is $\tilde M^3= \mathbb{R}^2\times \mathbb{R}_+$. The lift of $\nabla$ to $\tilde M^3$  is the Levi-Civita connection of the lift of $g$, which we denote by $\tilde g$. The reducibility of the holonomy group produces two (mutually orthogonal totally geodesic)  foliations on $\tilde M^3$. It is easy to see that
 \begin{enumerate} 
 \item 
 The leaf of one of these foliations is isometric to $\mathbb{R}$; the induced metric on it is complete and flat. 
 \item  
 The leaf of another foliation is $\mathbb{R}\times \mathbb{R}_+ $, and the induced metric on it is neither complete nor flat. 
 \item 
 The universal cover $\tilde M^3$ is isometric to the direct product of a leaf of the first foliation and a leaf of the second foliation (with the induces metrics).
 \end{enumerate}

We can show that the properties $(1,2,3)$ are necessary for any (real-analytic) counterexample to the Belgun-Moroianu conjecture; more precisely the following statement is true.
\emph{Suppose a closed connected manifold $M$ carries a reducible, locally metric, nonflat, real-analytic connection $\nabla$ which  preserves a conformal structure and which is not the Levi-Civita connection of any Riemannian metric. Then the universal cover  $\tilde M$ is isometric to the direct product of $(\mathbb{R}^k, g_{\mathrm{standard}})\times (N,h)$ such that the Riemannian manifold $(N, h)$ is neither complete, nor flat, and such that the lift of $\nabla$ is the Levi-Civita connection of the product metric $g_{\mathrm{standard}} + h$}.

\vspace{1ex} 
This statement, and also our study of reducible locally metric connections on closed manifolds were motivated by the study of conformally Berwald  Finlser manifolds. We will prove this statement, explain the motivation and show its applications in Finsler geometry elsewhere.

\end{remark}

\end{document}